\documentclass[12pt,a4paper,leqno]{amsart}

\usepackage{anysize}
\usepackage[utf8]{inputenc}
\usepackage{amsmath}
\usepackage{amssymb}
\usepackage{array}

\marginsize{3cm}{3cm}{3.5cm}{3.5cm}

\newcommand{\hpg}[5]{{}_{#1}\mbox{\rm F}_{\!#2}\! \left(\left.{#3 \atop #4}\right| #5 \right)}

\begin{document}

\title{New proofs of Borwein-type algorithms for Pi}
\author{Jesús Guillera}
\address{Department of Mathematics, University of Zaragoza, 50009 Zaragoza, SPAIN}
\email{jguillera@gmail.com}
\date{}

\maketitle

\begin{abstract}
We use a method of translation to recover Borweins' quadratic and quartic iterations. Then, by using the WZ-method, we obtain some initial values which lead to the limit $1/\pi$. We will not use the modular theory nor either the Gauss' formula that we used in \cite{Gui}. Our proofs are short and self-contained.
\end{abstract}

\section{Introduction}
The complete elliptic integral $K(x)$ is defined as follows:
\[
K(x)=\int_{0}^{\frac{\pi}{2}} \frac{dt}{\sqrt{1-x^2 \sin^2 t}},
\]
and its expansion in power series is given by
\[
\frac{2}{\pi}  K(x)=\hpg21{ \frac12, \, \frac12}{1}{x^2} = \sum_{k=0}^{\infty} \frac{\left(\frac12\right)_k^2}{(1)_k^2} x^{2k},
\]
where $(a)_k=a(a+1)\cdots(a+k-1)$ is the rising factorial, also known as Pochhammer symbol. Following the Indian mathematical genius Srinivasa Ramanujan, an identity of the form 
\[
\frac{K(\sqrt{1-v^2})}{K(v)}=n \frac{K(\sqrt{1-u^2})}{K(u)}, 
\]
where $n$ is any positive rational number, induces an algebraic transformation of degree $n$, which is of the following form \cite{ramanujan} or \cite[Ch. 29]{BeBo}:
\[
K(u)=m(u,v) \, K(v), \quad P(u,v)=0,
\]
where $m(u, v)$ (the multiplier) and $P(u, v)$ (the modular equation) are algebraic functions in $u$ and $v$. Many modular equations stated by Ramanujan in the variables $\alpha=u^2$ and $\beta=v^2$ are proved by Berndt in \cite[Ch.19]{Be3}. In  addition, Ramanujan proved that if $u_0$ is a solution of the equation $P(u, \sqrt{1-u^2})=0$, then there exist two algebraic numbers $C_1$ and $C_2$ such that 
\begin{equation}\label{rama-K2}
\left( C_1 \frac{K^2}{\pi^2}+C_2 \frac{K K'}{\pi^2} \right) (u_0)=\frac{1}{\pi},
\end{equation}
Observe that as $K$ is a common factor, we can write (\ref{rama-K2}) in the following way 
\[
\sum_{k=0}^{\infty} \frac{\left(\frac12\right)_k^2}{(1)_k^2} u_0^{2k} \sum_{k=0}^{\infty} \frac{\left(\frac12\right)_k^2}{(1)_k^2} u_0^{2k} (A+B k) = \frac{1}{\pi},
\]
where $A$ and $B$ are algebraic numbers related to $u_0$, $C_1$ and $C_2$. Instead, in his proof Ramanujan relates (\ref{rama-K2}) to a unique $_3F_2$ hypergeometric sum of the following form:
\[
\sum_{k=0}^{\infty} \frac{\left(\frac12\right)_k^3}{(1)_k^3} \left[4 u_0^2 (1-u_0^2)\right]^k (a+bk) = \frac{1}{\pi},
\]
where $a$, $b$ are algebraic numbers. One example given by Ramanujan is
\[
\sum_{k=0}^{\infty} \frac{\left(\frac12\right)_k^3}{(1)_k^3} \frac{5+42k}{64^k}=\frac{16}{\pi}.
\]
In the middle eighties, the brothers Jonathan and Peter Borwein observed that the modular equations given by Ramanujan could be connected to the number $\pi$ via the {\it elliptic alpha function}, to construct extraordinary rapid algorithms to calculate $\pi$ \cite[p.170]{Bo}. The point of view of the Borweins' is explained in \cite{Bo2}, and their proofs require a good knowledge of the elliptic modular functions. In this paper we show a different point of view, and although the idea of our method looks very general, we only focus in a quadratic transformation, which will allow us to recover the quadratic and quartic algorithms for a few sets of initial values (got with the WZ-method). The main point is that we avoid completely the use of the modular stuff, and also the Gauss's formula used in \cite{Gui}. 

\section{Transformations and algorithms}
Applying Zudilin's translation method \cite{Zu} to a quadratic transformation we will derive Borweins' quadratic iteration.

\subsection{A quadratic transformation}
The following identity: 
\begin{equation}\label{quad-trans}
\sum_{k=0}^{\infty} \frac{\left(\frac12\right)_k^2}{(1)_k^2}x^{2k}=(1+t) \sum_{k=0}^{\infty} \frac{\left(\frac12\right)_k^2}{(1)_k^2}t^{2k}.
\end{equation}
where
\begin{equation}\label{t-x}
t = \frac{1-\sqrt{1-x^2}}{1+\sqrt{1-x^2}},
\end{equation}
and $0<x<1$ is a known quadratic transformation \cite[Thm 1.2 (b)]{Bo}. Observe that the transformation is quadratic because $t \sim 4^{-1} \cdot x^2$. To prove (\ref{quad-trans})  just expand in powers of $x$, observe that the coefficients up to $x^2$ are equal and  check that both sides of (\ref{quad-trans}) satisfy the differential equation:
\[
\left( (x^2-1) \vartheta_x^2  + 2 x^2 \vartheta_x + x^2 \right) f = 0, \qquad \vartheta_x = x \frac{d}{dx},
\]
taking into account that the recurrence satisfied by the coefficient $D_k=(1/2)_k^2/(1)_k^2$ of $x^{2k}$ and $t^{2k}$ is:
\[
D_0=1, \qquad D_{k+1}=\frac{(1+2k)^2}{4(1+k)^2} D_k,
\]
the following relations which come from  (\ref{t-x}):
\[
\sqrt{1-x^2}=\frac{1-t}{1+t}, \qquad x^2=\frac{4t}{(1+t)^2}, \qquad x \frac{dt}{dx}=\frac{x^2(1+t)^3}{2-2t}=\frac{2t (1+t)}{1-t},
\]
and 
\[
\vartheta_x=\frac{x}{t} \, \frac{dt}{dx} \, \vartheta_t = \frac{2+2t}{1-t} \vartheta_t, \qquad 
\vartheta_x^2 = \frac{4(1+t)^2}{(1-t)^2}\vartheta_t^2 + \frac{8t(1+t)}{(1-t)^3} \vartheta_t.
\]

\subsection{A quadratic algorithm}\label{sec-quad-algo}
We consider the following operator:
\[
a+\frac{b}{2} \, \vartheta_x = a+b \, \frac{1+t}{1-t} \, \vartheta_t,
\]
and apply it to both sides of (\ref{quad-trans}):
\[
\left( a+\frac{b}{2} \, \vartheta_x \right) \left\{ \sum_{k=0}^{\infty} \frac{\left(\frac12\right)_k^2}{(1)_k^2}x^{2k} \right\} =
\left( a+b \, \frac{1+t}{1-t} \, \vartheta_t \right) \left\{ (1+t) \sum_{k=0}^{\infty} \frac{\left(\frac12\right)_k^2}{(1)_k^2}t^{2k} \right\}.
\]
We have
\begin{equation}\label{comp-trans}
\sum_{k=0}^{\infty} \frac{\left(\frac12\right)_k^2}{(1)_k^2} (a+bk) x^{2k}=
\sum_{k=0}^{\infty} \frac{\left(\frac12\right)_k^2}{(1)_k^2} 
\left( a(1+t)+b\frac{t(1+t)}{1-t} + 2b\frac{(1+t)^2}{1-t} k \right) t^{2k}
\end{equation}
Multiplying (\ref{quad-trans}) and (\ref{comp-trans}), we obtain
\begin{multline}\label{final-trans}
\sum_{k=0}^{\infty} \frac{\left(\frac12\right)_k^2}{(1)_k^2}x^{2k} \sum_{k=0}^{\infty} \frac{\left(\frac12\right)_k^2}{(1)_k^2} (a+bk) x^{2k} \\ =
\sum_{k=0}^{\infty} \frac{\left(\frac12\right)_k^2}{(1)_k^2}t^{2k} \sum_{k=0}^{\infty} \frac{\left(\frac12\right)_k^2}{(1)_k^2} 
\left( a(1+t)^2+b\frac{t(1+t)^2}{1-t} + 2b\frac{(1+t)^3}{1-t} k \right) t^{2k}
\end{multline}
Let
\[
A_n = \sum_{k=0}^{\infty} \frac{\left(\frac12\right)_k^2}{(1)_k^2} d_n^{2k}\sum_{k=0}^{\infty} \frac{\left(\frac12\right)_k^2}{(1)_k^2} (a_n+b_nk) d_n^{2k}.
\]
We define the sequences $d_n$, $a_n$, $b_n$, $A_n$ of initial values $d_0$, $a_0$, $b_0$, $A_0$, implicitly by means of $A_n=A_{n+1}$, where $A_{n+1}$ comes from $A_n$ in the same way than the right hand side of (\ref{final-trans}) comes from its left side. Hence we have the following recurrences:
\begin{align}
d_{n+1} &= \frac{1-\sqrt{1-d_n^2}}{1+\sqrt{1-d_n^2}}, \qquad b_{n+1} = 2 b_n \, \frac{(1+d_{n+1})^3}{1-d_{n+1}}, \label{ite-b}  \\
a_{n+1} &= a_n (1+d_{n+1})^2 + \frac{b_{n+1} d_{n+1}}{2(1+d_{n+1})} \label{ite-a}.
\end{align}
As $A_n$ is a constant sequence, it implies that $\lim A_n=A_0$, and as $d_n^{2k} \to 0$ and $b_n d_n^{2k} \to 0$ when $n \to \infty$ and $k \neq 0$, we deduce that $\lim A_n = \lim a_n$. Hence $\lim a_n = A_0$. 

\subsection*{Remark}
\rm 
Instead of dealing with algebraic transformations of products of two $_2F_1$ hypergeometric sums, we can deal with transformations of single $_3F_2$ hypergeometric series. In fact, there is an equivalence. For example, from the known transformation
\[
\left(\sum_{k=0}^{\infty} \frac{(\frac12)_k^2}{(1)_k^2}x^{2k} \right)^2 = \sum_{k=0}^{\infty} \frac{(\frac12)_k^3}{(1)_k^3}\left[ 4x^2(1-x^2) \right]^k,
\]
we can, via translation, rewrite our products of two $_2F_1$ series as unique $_3F_2$ series. When the sum is equal to $1/\pi$ these single sums constitute one of the families of Ramanujan-type series for $1/\pi$. 

\section{Other algorithms}
Here we give some related algorithms:

\subsection{A simplified version}
From (\ref{ite-b}), we have
\[
\frac{b_{n+1}}{1-d_{n+1}^2} = 2 b_n \frac{(1+d_{n+1})^2}{(1-d_{n+1})^2} = 2 \frac{b_n}{1-d_n^2}.
\]
Then, writing $c_n=b_n / (1-d_n^2)$, we have $c_{n+1}=2c_n$, which imply $c_n = c_0 \cdot 2^n$, and we get a simplified version of the algorithm of subsection \ref{sec-quad-algo}. It is
\begin{align}
d_{n+1} &= \frac{1-\sqrt{1-d_n^2}}{1+\sqrt{1-d_n^2}}, \label{ite2-d} \\
a_{n+1} &= a_n (1+d_{n+1})^2 + c_0 \, 2^n d_{n+1} (1-d_{n+1}) \label{ite2-a},
\end{align}
with $\lim a_n = A_0$. 

\subsection{An even more simplified version: Borweins' quadratic algorithm}
As $c_0 2^n d_{n+1} \sim a_{n+1}-a_n$, we deduce that $2^n d_{n+1}$ tends to $0$ as $n \to \infty$. Then, writing $a_n=r_n-c_0 \,  2^{n-1} \, d_n^2$, we arrive to the Borweins' quadratic iteration
\begin{align}
d_{n+1} &= \frac{1-\sqrt{1-d_n^2}}{1+\sqrt{1-d_n^2}}, \label{ite3-d} \\
r_{n+1} &= r_n (1+d_{n+1})^2 - c_0 \, 2^n d_{n+1} \label{ite3-r},
\end{align}
with $\lim r_n = \lim a_n =A_0$. 

\subsection{Borweins' quartic algorithm}
Making $s_n=\sqrt{d_{2n}}$, $t_n=r_{2n}$, we can arrive to the following Borweins' quartic algorithm (see \cite{Gui} for the steps):
\begin{align*}
s_0 &=\sqrt{d_0}, \quad t_0=r_0, \\
s_{n+1} &= \frac{1-\sqrt[4]{1-s_n^4}}{1+\sqrt[4]{1-s_n^4}}, \\
t_{n+1} &=(1+s_{n+1})^4 t_n - c_0 2^{2n+1} s_{n+1} (1+s_{n+1}+s_{n+1}^2), 
\end{align*}
with $\lim t_n = A_0$. Observe that it is a quartic algorithm because $s_{n+1} \sim 8^{-1} \cdot s_n^4$.

\subsection{Initial values}
The identities of this subsection provide initial values in order that the iterations tend to $A_0=1/\pi$. Below each identity we show the corresponding initial values in the algorithms. We will prove these evaluations in next subsection. 

\subsubsection*{Identity 1}
\[
\sum_{k=0}^{\infty} \frac{(\frac12)_k^2}{(1)_k^2} (3-2\sqrt 2)^{2k}
\sum_{k=0}^{\infty} \frac{(\frac12)_k^2}{(1)_k^2} (3-2\sqrt 2)^{2k} \left[ 16(3\sqrt 2 - 4) k + 4(5\sqrt 2 -7) \right]= \frac{1}{\pi}, 
\]
and
\[
d_0 =s_0^2=3-2 \sqrt 2, \quad b_0=48\sqrt 2-64, \quad a_0=20\sqrt 2-28, \quad c_0=4, \quad r_0=t_0=6-4 \sqrt 2.
\]

\subsubsection*{Identity 2}
\[
\sum_{k=0}^{\infty} \frac{(\frac12)_k^2}{(1)_k^2} \, \frac{1}{2^k}
\sum_{k=0}^{\infty} \frac{(\frac12)_k^2}{(1)_k^2} \, \frac{k}{2^k} = \frac{1}{\pi},
\]
and
\[
d_0=s_0^2=\frac{1}{\sqrt 2}, \quad b_0=1, \quad a_0=0, \quad c_0=2, \quad r_0=t_0=\frac12.
\]

\subsubsection*{Identity 3}
\[
\sum_{k=0}^{\infty} \frac{(\frac12)_k^2}{(1)_k^2} (\sqrt 2-1)^{2k}
\sum_{k=0}^{\infty} \frac{(\frac12)_k^2}{(1)_k^2} (\sqrt 2-1)^{2k} \left[ (8-4\sqrt 2) k + (3-2\sqrt 2) \right]= \frac{1}{\pi}, 
\]
and
\[
d_0 =s_0^2=\sqrt 2 - 1, \quad b_0=8-4 \sqrt 2, \quad a_0=3-2\sqrt 2, \quad c_0=2 \sqrt 2, \quad r_0 =t_0=\sqrt 2 - 1.
\]

\subsubsection*{Identity 4}
\[
\sum_{k=0}^{\infty} \frac{(\frac12)_k^2}{(1)_k^2} \left(\frac{2-\sqrt 3}{4} \right)^k
\sum_{k=0}^{\infty} \frac{(\frac12)_k^2}{(1)_k^2} \left(\frac{2-\sqrt 3}{4} \right)^k \left[ \left(\frac32 + \sqrt 3 \right) k + \frac14 \right] = \frac{1}{\pi},
\]
and
\[
d_0 =s_0^2=\frac{\sqrt 6 - \sqrt 2}{4}, \quad b_0=\frac32+\sqrt 3, \quad a_0=\frac14, \quad c_0=2 \sqrt 3, \quad r_0 =t_0 = \frac{\sqrt 3 - 1}{2}. 
\]

\subsubsection*{Identity 5}
\[
\sum_{k=0}^{\infty} \frac{(\frac12)_k^2}{(1)_k^2} \left(2\sqrt 2 - 2 \right)^k
\sum_{k=0}^{\infty} \frac{(\frac12)_k^2}{(1)_k^2} \left(2\sqrt 2 -2  \right)^k \left[ (3\sqrt 2-4) k + \sqrt 2 - \frac32 \right] = \frac{1}{\pi},
\]
and
\[
d_0 =s_0^2=\sqrt{2\sqrt 2-2}, \quad b_0=3\sqrt 2 - 4, \quad a_0=\sqrt 2 - \frac32, \quad c_0=\sqrt 2, \quad r_0 = t_0=\frac12.
\]

\subsection{Elementary proofs of identities 1 to 5}
In this section we prove Identities 1 to 5. We avoid the modular theory and use only the WZ-method together with Carlson's Theorem \cite{ekhad}. In the right side of the formulas we use the general definition of the Pochhammer symbol, namely $(v)_k=\Gamma(v+k)/\Gamma(v)$.

\subsubsection*{WZ-proof of Identity 1}
The following identity: 
\begin{equation}\label{2F1-1}
\sum_{n=0}^{\infty} \frac{(\frac12-4k)_n (\frac12-2k)_n}{(1+2k)_n(1)_n} (3-2\sqrt 2)^{2n} = C_1 (3\sqrt 2-4)^{2k} \frac{2^{9k}}{3^{3k}} \, \frac{(1)_k \left(\frac12\right)_k}{\left(\frac{7}{12}\right)_k (\frac{11}{12})_k},
\end{equation}
valid for any complex number $k$ is WZ-provable \cite{PWZ}, due to the fact that it has a free parameter, and we can determine the constant $C_1$ by taking $k=1/4$. To discover other formulas like this, we need that the function $H(n,k)$ inside the sum is such that the output of the Maple code
\begin{verbatim}
                with(SumTools[Hypergeometric]);
                degree(Zeilberger(H(n,k),k,n,K)[1],K);
\end{verbatim}
is $1$. In this case we can easily guess the function at the right hand side of the equal, which is a function depending only in $k$. Many terminating versions of formulas of this kind have been discovered and proved automatically by a program written by D. Zeilberger \cite{ekhad2}. It is good for our purposes that we found the following complement of (\ref{2F1-1}):
\begin{align}
\sum_{n=0}^{\infty} \frac{(\frac12-4k)_n (\frac12-2k)_n}{(1+2k)_n(1)_n} (3-2\sqrt 2)^{2n} & \left[ 16(3\sqrt 2 -4) n + 16(13-9 \sqrt 2) k + 4(5\sqrt 2 -7) \right] \nonumber \\ &= C_2 (3\sqrt 2-4)^{2k} \frac{2^{9k}}{3^{3k}} \, \frac{(1)_k \left(\frac12\right)_k}{\left(\frac{1}{12}\right)_k (\frac{5}{12})_k}, \label{2F1-2}
\end{align}
which is WZ-provable as well. To discover this complement we have used the following Maple code: 
\begin{verbatim}
       with(SumTools[Hypergeometric]);
       coK2:=coeff(Zeilberger(H(n,k)*(n+b*k+c),k,n,K)[1],K,2);
       coes:=coeffs(coK2,k); solve({coes},{b,c});
\end{verbatim}
where $H(n,k)$ is the function inside the sum in (\ref{2F1-1}) written in the notation of Maple. Again, we can determine $C_2$ taking $k=1/4$, and  we have that $C_1 C_2=1/\pi$. Multiplying the formulas (\ref{2F1-1}) and (\ref{2F1-2}), taking $k=0$, and finally replacing $n$ with $k$ to agree with the notation of the other subsections, we arrive at Identity 1.
\vskip 0.25cm
Below we give  couples of identities which prove Identities 2 to 5.

\subsubsection*{WZ-proof of Identity 2}
\begin{align*}
& \sum_{n=0}^{\infty} \frac{(\frac12-2k)_n (\frac12+2k)_n}{(1+4k)_n(1)_n} \left(\frac12\right)^n = C_1 \left( \frac{16}{27} \right)^k \frac{(1)_k \left(\frac12\right)_k}{\left(\frac{7}{12}\right)_k (\frac{11}{12})_k}, \hskip 5cm \\
& \sum_{n=0}^{\infty} \frac{(\frac12-2k)_n (\frac12+2k)_n}{(1+4k)_n(1)_n} \left(\frac12\right)^n (n+4k)  = C_2 \left( \frac{16}{27} \right)^k \frac{(1)_k \left(\frac12\right)_k}{\left(\frac{1}{12}\right)_k (\frac{5}{12})_k}. 
\end{align*}

\subsubsection*{WZ-proof of Identity 3}
\begin{align*}
& \sum_{n=0}^{\infty} \frac{(\frac12-2k)_n (\frac12)_n}{(1+2k)_n(1)_n} (\sqrt 2-1)^{2n} = C_1 \, \frac{(1)_k \left(\frac12\right)_k}{\left(\frac{5}{8}\right)_k (\frac{7}{8})_k}, \\
& \sum_{n=0}^{\infty} \frac{(\frac12-2k)_n (\frac12)_n}{(1+2k)_n(1)_n} (\sqrt 2-1)^{2n} \left[ (8-4\sqrt 2)n + (8\sqrt 2-8)k + (3-2\sqrt 2) \right] = C_2 \, \frac{(1)_k \left(\frac12\right)_k}{\left(\frac{1}{8}\right)_k (\frac{3}{8})_k}.
\end{align*}

\subsubsection*{WZ-proof of Identity 4}
\begin{align*}
& \sum_{n=0}^{\infty} \frac{(\frac12-k)_n (\frac12+3k)_n}{(1+k)_n(1)_n} \left( \frac{2-\sqrt 3}{4}\right)^n = C_1 \, \left( \frac{4}{3\sqrt 3} \right)^k \, \frac{(1)_k}{\left(\frac56\right)_k}, \\
& \sum_{n=0}^{\infty} \frac{(\frac12-k)_n (\frac12+3k)_n}{(1+k)_n(1)_n} \left( \frac{2-\sqrt 3}{4}\right)^n \left[ \left(\frac32+\sqrt 3\right)n + \frac32 k + \frac14 \right] = C_2 \, \left( \frac{4}{3\sqrt 3} \right)^k \, \frac{(1)_k}{\left(\frac16 \right)_k}.
\end{align*}

\subsubsection*{WZ-proof of Identity 5}
\begin{align*}
& \sum_{n=0}^{\infty} \frac{(\frac12+4k)_n (\frac12+2k)_n}{(1+k)_n(1)_n} (2 \sqrt 2-2)^n = C_1 \, (3+2 \sqrt 2)^{2k} \, \frac{(1)_k \left(\frac12\right)_k}{\left(\frac58\right)_k\left(\frac78\right)_k}, \hskip 3.25cm \\
& \sum_{n=0}^{\infty} \frac{(\frac12+4k)_n (\frac12+2k)_n}{(1+k)_n(1)_n} (2 \sqrt 2-2)^n \left[ (3\sqrt 2 -4)n+(8\sqrt 2 -12)k+\left(\sqrt 2 - \frac32\right) \right] \\ & \hskip 6.25cm = C_2 \, (3+2 \sqrt 2)^{2k} \, \frac{(1)_k \left(\frac12\right)_k}{\left(\frac18\right)_k \left(\frac38\right)_k}. 
\end{align*}

Although with the WZ-method we have only obtained a few sets of initial values, it supposes no disadvantage with respect to the modular theory because the interesting fact is to have one set of simple initial values.

\end{document}